\documentclass[reqno,12pt]{amsart}

\usepackage{amsmath}
\usepackage{amssymb}
\usepackage{amsthm}
\usepackage[english]{babel}

\usepackage{fancyhdr}  
\newtheorem{lemma}{Lemma}

\begin{document}

\title {Bound for the maximal probability in the
Littlewood--Offord problem}

\author{Andrei Yu. Zaitsev}


\footnotetext[1]{The research is supported by grants RFBR
13-01-00256 and NSh-2504.2014.1 and by the Program of Fundamental
Researches of Russian Academy of Sciences "Modern Problems of
Fundamental Mathematics".}

\email{zaitsev@pdmi.ras.ru}
\address{St.~Petersburg Department of Steklov Mathematical Institute
\newline\indent
Fontanka 27, St.~Petersburg 191023, Russia\newline\indent and
St.~Petersburg State University}

\begin{abstract}{The paper deals with studying a connection of the
Littlewood--Offord problem with estimating the concentration
functions of some symmetric infinitely divisible distributions. It
is shown that the values at zero of the concentration functions of
weighted sums of i.i.d. random variables may be estimated by the
values at zero of the concentration functions of symmetric
infinitely divisible distributions with the L\'evy spectral
measures which are  multiples of the sum of delta-measures at
$\pm$weights involved in constructing the weighted
sums.}\end{abstract}

\keywords {concentration functions, inequalities, the
Littlewood--Offord problem, sums of independent random variables}

\maketitle

The
 concentration function of a $\mathbf{R}^d$-dimensional random
vector $Y$ with distribution $F=\mathcal L(Y)$ is defined by the
equality
$$
Q(F,\tau)=\sup_{x\in\mathbf{R}^d}\mathbf{P}(Y\in x+ \tau B), \quad
\tau\ge0, \nonumber
$$
where $B=\{x\in\mathbf{R}^d:\|x\|\leq 1/2\}$ is the centered ball
of the Euclidean space $\mathbf{R}^d$ of radius 1/2. In
particular,
$$
Q(F,0)=\sup_{x\in\mathbf{R}^d}\mathbf{P}(Y= x). \nonumber
$$

Let $X,X_1,\ldots,X_n$ be independent identically distributed
 random variables. Let $a=(a_1,\ldots,a_n)$, where
$a_k=(a_{k1},\ldots,a_{kd})\in \mathbf{R}^d$, $k=1,\ldots, n$.
Starting with seminal papers of Littlewood and Offord
\cite{Lit:Off:1943} and Erd\"os~\cite{Erd}, the behavior of the
concentration functions of the weighted sums
$S_a=\sum\limits_{k=1}^{n} X_k a_k$ is studied intensively. In the
sequel, let $F_a$ denote the distribution of the sum $S_a$. The
first version of the Littlewood--Offord problem was related to the
estimation of $Q(F_a, 0)$ in the case where $X$ has the symmetric
Bernoulli distribution $\mathbf{P}(X= 1)=\mathbf{P}(X= -1)=1/2$.

In the last ten years, refined concentration results for the
weighted sums $S_a$ play an important role in the study of
singular values of random matrices (see, for instance, Nguyen and
Vu~\cite{Nguyen:Vu:2011, Nguyen and Vu13}, Rudelson and Vershynin
\cite{Rud:Ver:2008, Rud:Ver:2009}, Tao and Vu
\cite{Tao:Vu:2009:Ann, Tao:Vu:2009:Bull},
Vershynin~\cite{Ver:2011}, Later, somewhat different bounds for
the concentration functions in the Littlewood--Offord problem were
obtained by Eliseeva, G\"otze and Zaitsev
\cite{Eliseeva}--\cite{EZ}.
 The aforementioned results reflect the dependence
of the bounds on the arithmetic structure of coefficients~$a_k$
under various conditions on the vector $a\in {(\mathbf{R}^d)}^n$
and on the distribution~$\mathcal L(X)$.

Below, in some formulas, the quantity
$$p(v)=
G\big\{\{x\in\mathbf{R}:|x| > v\}\big\},\quad v\ge0,
$$ will appear,
where $G$ is the symmetrized distribution $G={\mathcal
L}(X_1-X_2)$. Writing ${A\ll _{d} B}$ means that $|A|\leq c(d)\,
B$, where $c(d)>0$ depends on $d$ only. The inner product of
vectors $x,y\in\mathbf{R}^d$ is written as $\left\langle
x,y\right\rangle$. Denote by $\frak F_d$ the set of all
$d$-dimensional probability distributions and by $\frak D_d$ the
set of all infinitely divisible distribution from $\frak F_d$. We
denote by $E_y$ the distribution concentrated at a point $y\in
\mathbf{R}^{d}$.
 Below $\widehat F(t)$, $t\in\mathbf
R^d$,  is the characteristic function of~$F\in\frak F_d$.

Products and powers of measures will be understood in the sense of
convolution. For $D\in\frak D_d$, $F\in\frak F_d$, and
$\lambda\ge0$, by $D^\lambda$ we denote infinitely divisible
distribution with the characteristic function $\widehat
D^\lambda(t)$, and by ${e}(\lambda\, F)$ infinitely divisible
compound Poisson distribution with the characteristic function
\,$\exp \big(\lambda\,(\widehat{F}(t)-1)\big)$, and with the
L\'evy spectral measure  $\lambda\, F$. It is easy to see that
$$
{e}(\lambda\,F)=e^{-\lambda}\sum_{s=0}^{\infty}\frac{\lambda^{s}F^{s}}{s!}.
$$
Here $F^0$ is the degenerate distribution $E_0$ concentrated at
the origin $0\in \mathbf{R}^{d}$. Thus,
\begin{equation}\label{98t}
{e}(\lambda\,F)=e^{-\lambda}E_0+
e^{-\lambda}\sum_{s=1}^{\infty}\frac{\lambda^{s}F^{s}}{s!}.
\end{equation}

We need some simplest properties of concentration functions. First
we note that for any distribution $F\in\frak F_d$
\begin{equation}\label{8rt}
\lim_{\tau\to0} Q(F, \tau) = Q(F, 0)
\end{equation}
(see, e.g.,  \cite[p. 14]{Hen:Teod:1980}). For any $U,V\in\frak
F_d$, we have
\begin{equation}\label{8t}
 Q(\,U \,V, \tau) \le Q(U, \tau),\quad\hbox{for all }\tau\ge0.
\end{equation}

The following Lemma \ref{lm3} follows directly from \eqref{98t}
and \eqref{8t}.

\begin{lemma}\label{lm3} Let $\tau,\lambda\geq0$,
 $F, U\in\frak F_d$ and $D={e}(\lambda\,F)\in\frak D_d$. Then
\begin{equation}
0\leq Q(U, \tau)-Q(U D, \tau)\le 1-e^{-\lambda}.
 \label{4s}
\end{equation}
\end{lemma}
 \medskip

 Introduce the
distribution $H$, with the characteristic function
$$\widehat H(t)
=\exp\Big(-\frac{1}2\;\sum_{k=1}^{n}\big(1-\cos\langle \,
t,a_k\rangle\big)\Big),\quad t\in \mathbf{R}^d. $$ Note that
$H^p={e}\big(\frac {\,n\,p\,}2\,M\big)$ is a symmetric infinitely
divisible distribution with the L\'evy spectral measure $\frac
{\,n\,p\,}2\,M$, where $$M= \frac
1{2\,n}\sum_{k=1}^{n}\big(E_{a_k}+E_{-a_k}\big)\in\frak F_d.$$ M.
A. Lifshits drew the author's attention to the fact that the
distribution $H^p$ can be represented as the distribution of
weighted sum $S_a$ for the same set of weights
$a=(a_1,\ldots,a_n)$ which is involved in the the original
problem. But the common distribution of the random variables
$X,X_1,\ldots,X_n$, for fixed $p$, has a special form ${\mathcal
L}(X)={e}\big(\frac {\,p\,}4\,(E_{1}+E_{-1}\big)\big)$.

The following Theorem 1 is contained in the recently published
paper \cite{EZ}, see also \cite{EGZ1}. It connects the
Littlewood--Offord problem with general bounds for the
concentration functions, in particular, with the results of Arak
contained in \cite{Arak:1980}--\cite{Arak:Zaitsev:1986} (see
\cite{EGZ1}).
 \medskip

\noindent {\bf Theorem 1.} {\it For any $\tau,u>0$, the inequality
$$
Q(F_a, \tau) \ll_d Q(H^{p(\tau/u)}, u) $$ holds.}
 \medskip

In a recent paper of Eliseeva and Zaitsev \cite{EZ}, a more
general statement than Theorem 1 is obtained. It gives useful
bounds if $p(\tau/u)$ is small, even if $p(\tau/u)=0$.

Theorem 1 has been proved for $\tau, u> 0$. The naturally arising
question is to find an analogue of Theorem 1 for $\tau = 0$. The
answer to this question is given by the following Theorem 2.
 \medskip

\noindent {\bf Theorem 2.} {\it The inequality
\begin{equation}\label{978t}
Q(F_a, 0) \ll_d Q(H^{p(0)}, 0)   \end{equation} holds.} \medskip

Let $A=\big\{x\in\mathbf R:\mathbf P(X=x)>0\big\}$. It is clear
that the set $ A $ is no more than countable. It is easy to verify
that
$$p(0)=1-\sum_{x\in A}\big(\mathbf P(X=x)\big)^2.$$
According to \eqref{8t}, $Q(H^{p(0)}, 0)$ is a non-increasing
function of $p(0)$. Therefore, the right-hand side of inequality
\eqref{978t} is minimal for $ p (0) = 1 $. But in this case the
left-hand side of \eqref{978t} vanishes. Indeed, then the
distribution ${\mathcal L}(X)$ has no atoms, $A=\varnothing$,  and
$Q\big({\mathcal L}(X), 0\big)=0$. Through \eqref{8t}, this easily
implies the equality $Q(F_a, 0)=0$.

The advantage of Theorem 2 is that  $p (v)$ takes its maximal
value at $ v = 0 $. It distinguishes  Theorem 2 from Theorem 1,
which gives for $\tau=u>0$ the estimate $Q(F_a, \tau) \ll_d
Q(H^{p(1)}, \tau)$,  which, of course, implies that $Q(F_a, 0)
\ll_d Q(H^{p(1)}, 0)$ (see \eqref{8rt}). But $ p (0) $ can be
substantially greater than $ p (1) $ and $Q(H^{p(0)}, 0)$ may be
substantially less than $Q(H^{p(1)}, 0)$. Using Theorem~2 instead
of Theorem 1 gives a possibility to replace $p(1)$ by $p (0)$ in
Theorems 5 and 6 of \cite{EGZ1} in a particular case, where the
parameters $\tau_j$, $j=1,\dots,d$, involved in the formulations
of these theorems, are all zero.

It is interesting that, in spite of the above, we deduce Theorem 2
still from Theorem 1 with the help of somewhat more delicate
passing to the limit.
\medskip

\noindent{\bf Proof of Theorem 2.} By Theorem 1,
$$
Q(F_a, \tau) \ll_d Q(H^{p(\varepsilon)},
\tau/\varepsilon),\quad\hbox{for all }\tau,\varepsilon>0. $$
Letting $\tau$ to zero, and using \eqref{8rt}, we get
\begin{equation}\label{8t4}  Q(F_a,0) \ll_d Q(H^{p(\varepsilon)},
0),\quad\hbox{for all }\varepsilon>0. \end{equation}

It is evident that $p(\varepsilon)\le p(0)$ and $H^{p(0)}=
H^{p(\varepsilon)}H^{p(0)-p(\varepsilon)}$. Using
inequality~\eqref{4s}, we verify the validity of the relation
$$0\leq Q(H^{p(\varepsilon)}, 0)- Q(H^{p(0)},0)\le 1-\exp\Big(-\frac
1{\,2\,}\,n\,\big(p(0)-p(\varepsilon)\big)\Big) . $$ Since
$p(\varepsilon)\to p(0) $, this implies that
$$ Q(H^{p(\varepsilon)},
0)\to Q(H^{p(0)}, 0)\quad\hbox{as }\varepsilon\to0. $$ Letting
$\varepsilon$ to zero in the right-hand side of
inequality~\eqref{8t4}, we obtain the assertion of Theorem~2.
 \hfill$\square$
\medskip

\noindent{\bf Remark.} The weak convergence of probability
distributions does not imply, in general, the convergence of
values of the concentration functions at zero. For example, this
happens when continuous distributions converge to discrete ones.

\end{document}